\documentclass[12pt,oneside,a4paper,leqno]{article}

\usepackage{amsfonts,amssymb}
\usepackage{geometry}
\usepackage[centertags]{amsmath}
\usepackage[all]{xy}\xyoption{all}
\usepackage{graphicx}
\usepackage{amsthm}
\newcounter{lemma}
\swapnumbers
\theoremstyle{plain}
\newtheorem{lemma}[equation]{Lemma}\numberwithin{lemma}{section}

\newtheorem{propo}[equation]{Proposition}

\theoremstyle{definition}

\newtheorem{remark}[equation]{Remark}
\theoremstyle{remark}

\numberwithin{equation}{section}

\makeatletter
\def\tagform@#1{\maketag@@@{\ignorespaces#1\unskip\@@italiccorr}}

\makeatother

\newdir{ >}{!/8pt/@{ }*@{>}}
\newdir{< }{!/-8pt/@{ }*@{<}}

\theoremstyle{plain}
\newtheorem{assumption}[equation]{Collision Assumption}
\usepackage{paralist}

\newcommand{\C}{\mathbb{C}}
\newcommand{\R}{\mathbb{R}}

\newcommand{\from}{\colon}
\newcommand{\Chi}{\mathcal{X}}

\newcommand{\action}{\mathcal{A}}
\newcommand{\nn}{\mathbf{n}}

\newcommand{\st}{\ \mathrm{:} \ }
\newcommand{\ze}{\mathbb{Z}}
\newcommand{\Usf}{U_{\mathrm{sf}}}

\begin{document}
\pagenumbering{arabic}

\title{%
Symmetric periodic orbits for  the $n$-body problem:
some preliminary results
}

\author{Davide L.~Ferrario}

\date{}
\maketitle

\begin{abstract}
We show the existence of some infinite families of
periodic solutions of the planar Newtonian $n$-body problem --with
positive masses-- which are symmetric with
respect to suitable actions of finite groups
(under a strong--force assumption or only numerically).
The method is by minimizing a discretization of the action functional 
under symmetry constraints.
\end{abstract}




\section{Introduction}
\label{sec:intro}
Exact periodic solutions of the Newtonian $n$-body problem
have long been a topic of interest, and not only among the specialists
of celestial mechanics.
The recently found ``eight'' choreography of Montgomery and Chenciner
\cite{monchen}
has been the starting point of the numerical discovery
by Sim\'o and others of several 
interesting periodic orbits of the same kind,  
which are symmetric under the action
of a finite group \cite{chenven,chenciner}.
In this note we simply define suitable actions
of some finite groups on the configuration space
and infer the existence of periodic 
orbits 
by  the well-known variational approach. For details and 
further reading the reader is referred to 
\cite{ambbad,argate,mate,setate,mate95,
monchen,chenven,andrea,ambrosetticoti,ambrocoti}.
The most recent advances in this theory will appear
in Chenciner's papers for the ICM 2002 (see also \cite{chenciner}).
We need to consider the assumption that the local minima
of the Lagrangian action for the corresponding Bolza problem
is collision-free.
This has been proved to be true by Marshall (in a preprint) 
in case the symmetry group is cyclic and the potential is Newtonian;
it deserves some effort to be proved in general. But in any case
we can circumvent this problem by a strong--force perturbation of
the potential.
The paper is roughly structured as follows: we will prove rigorously 
(even if omitting full details, that can be find in the cited papers) 
some needed
results, and then use  numerical simulations to verify the 
collision assumption when needed and to compute the orbits.
We obtain some non--trivial periodic orbits in the following cases.

\begin{inparaenum}[\itshape (a)]
\item
$3$-body with all equal masses: the Montgomery--Chenciner choreography
(section \ref{subsec:threeeq}).
\item 
$3$-body with two equal masses:
infinitely many periodic orbits (section \ref{subsec:3_2eq} and 
figure \ref{fig:candidate}). 
Among these, there are 
the orbits with all the three masses equal (figure \ref{fig:3b}).
\item 
$4$-body with all equal masses:
the method only yields some numerical examples  (section \ref{subsec:4eq}), 
like the one in figure \ref{fig:nonhomographic}.
\item 
$4$-body with two pairs of equal masses:
an infinite family of solutions 
(section \ref{subsec:twopairs} and 
figure \ref{fig:rough}).
\item 
$4$-body with all equal masses and an additional symmetry:
an infinite family of periodic orbits (section \ref{subsec:4eqp},
figures \ref{fig:mac}, \ref{fig:caseq6}).
\item
$4$-body with two pairs of equal masses and an additional symmetry:
an infinite family (section \ref{subsec:22eqp},
figure \ref{fig:arch}).
\item 
Plane choreographies (``eight'' shape) with $n$ odd equal masses
(section \ref{sec:choreo}, figures \ref{fig:chor5}, \ref{fig:chor5bis}).
\end{inparaenum}

Many of these orbits are well-known (but in general only numerically),
other might not. This is an attempt ---most 
far from being ultimate--- to give a general method 
for proving the existence of such periodic orbits with variational
methods combined with a computer--assisted proof. More general
results (not confined to dihedral groups and without the collision
assumption) will come in subsequent papers.
Most of the ideas are borrowed from Susanna Terracini, who helped 
me to understand the problem and willingly shared  her knowledge
during several discussions on the topic. Sincere thanks
are also due to Andrea Venturelli, who generously gave his 
much appreciated comments.

\section{Preliminaries}
Let $n\geq 2$ be  the number of bodies in the 
Euclidean space $V=E^k \approx \R^k$ of dimension $k\geq 2$.  
Here we consider the $n$-body problem with masses $m_i$, $i=1,\dots,n$.
Following the notation of \cite{chenciner},
let $\Chi$ denote the subspace of $V^n$
of all points $x=(x_1,x_2,\dots,x_n)$ such that 
$\sum_{i=1}^n m_ix_i = 0$, i.e., with center of mass in the origin $O$.
Let $\Delta_{i,j}$ denote the collision set in $\Chi$ of the $i$-th
and the $j$-th particles,
and $\Delta = \cup_{i<j} \Delta_{i,j}$ the collision set
in $\Chi$.
Let $U\from \Chi \to \R$ be the potential function
\[
U(x)  = \sum_{i<j} \dfrac{1}{|x_i - x_j|}.
\]
We can consider also a deformation $\Usf$ of $U$ in a small 
neighborhood of the collision set, so that $U$ satisfies
the \emph{strong force} condition.
The kinetic energy, defined on the tangent bundle of $\Chi$,
is 
$K=\sum_{i=1}^n \dfrac{1}{2} m_i |\dot x_i|^2$, and the Lagrangian
is $L=K+U$. 
Moreover $I=\sum_{i} m_i |x_i|^2$ is the moment of inertia with 
respect to the center of mass.
Let $T^1\subset \R^2$ denote the unit circle and 
$\Lambda = H^1(T^1,\Chi)$ the Sobolev space 
of the $L^2$ loops $T^1 \to \Chi$ with $L^2$ derivative. 
Then, 
the critical points of the positive-defined action functional
$\action \from \Lambda \to \R\cup \{\infty\} $
are 
periodic orbits of the Newtonian $n$-body problem;
the action is defined by 
\begin{equation}\label{eq:action}
\action(x) = \int_{T^1} L(x(t),\dot x(t) ) dt
\end{equation}
for every loop $x=x(t) \in \Lambda$.
The action functional is called 
\emph{coercive} if it goes to infinity as 
the moment of inertia $I$ goes to infinity.
(i.e., if the $H^1$-norm of $x$ goes to infinity).

Let $G$ be a finite group, acting on a space $X$. The space $X$ is 
then called \emph{$G$-equivariant} space.
We recall some 
standard notation.
If $H\subset G$ is a subgroup of $G$, then 
$G_x = \{g\in G \st gx=x \}$ is termed the \emph{isotropy} of $x$,
or the  \emph{fixer} of $x$ in $G$.
The space $X^H\subset X$ consists of all the points 
$x\in X$ which are fixed by $H$, that is
$X^H=\{x\in X \st G_x \supset H \}$.
Given two $G$-equivariant spaces $X$ and $Y$,
an \emph{equivariant map} $f\from X \to Y$ is a map
with the property that $f(g\cdot x) = g\cdot f(x)$ 
for every $g\in G$ and every $x\in X$.
An equivariant map induces, by restriction to the spaces
$X^H$ fixed by subgroups $H\subset G$, maps
$f^H\from X^H \to Y^H$.

\section{Symmetry constraints}
We give here an introduction to symmetry
constraints which is slightly different than the well-known 
in the literature.
Let $G$ be a finite group, $\tau$ an orthogonal representation of 
$G$ on 
$T^1$ and  
$\rho$ an orthogonal representation on the Euclidean space
$V$. Furthermore, let $\sigma$ 
be a group homomorphism $\sigma\from G \to \Sigma_n$ 
from $G$ to the symmetric group on $n$ elements. Therefore
we let $G$ act on the space $V$,  on the time $T^1$ and on  the 
set of indexes of the masses $\mathbf{n}=\{1,2,\dots, n\}$.
Let $\R[\mathbf{n}]$ denote the 
$n$-dimensional real vector space generated by the $n$
elements of $\mathbf{n}$,
and 
$\R_0[\nn]$ the linear subspace 
of elements with coordinates $(\lambda_1,\lambda_2,\dots,\lambda_n)$
such that $\sum_{i=1}^n m_i\lambda_i = 0 $.
Given $\sigma\from G \to \Sigma_n$,  $G$ acts on 
$\R_0[\nn]$ by the action on $\nn$,
provided that $\sigma(g)(i) = j$ implies $m_i = m_j$
for every $g\in G$ and every $i=1,\dots, n$.
Two homomorphisms  $\sigma$ and $\sigma'$  yield
the same real representation on $\R_0[\nn]$
if they are conjugated by an element of $\Sigma_n$,
that is if one is obtained by permuting coordinates on the other.
With an abuse of notation we call $\sigma$ the representation
induced by the homomorphism $\sigma$. It is equivalent to 
the representation
given by $[N] - [1]$ (in the representation ring of $\Sigma_n$),
where $N$ is the natural representation and $1$ the trivial 
representation.
Now, since 
\[
\Chi \cong V\otimes_\R \R_0[\nn],
\]
given $\rho$ and $\sigma$ we have 
an orthogonal action of $G$ on $\Chi$.

Furthermore, since $G$ acts on $T^1$ 
and on $\Chi$, then there is the  standard diagonal  action
on the loop space $\Lambda$, defined by
$x(t) \mapsto g x( g^{-1} t)$. 
Let us note that 
the loops in $\Lambda$ fixed by $G$ are 
the equivariant maps $x\from T^1 \to \Chi$.
In this terminology, a \emph{symmetry constraint} is 
a such action of $G$ on $\Lambda$. Since the action
functional is invariant with respect to 
the $G$-action, 
we have a restricted action 
\begin{equation}
\label{eq:lambdag}
\action^G\from \Lambda^G \subset \Lambda \to \R,
\end{equation}
and 
the following proposition (Palais principle of symmetric
criticality -- see \cite{chenciner}).

\begin{lemma}\label{palais}
A critical point of $\action^G$ in $\Lambda^G$ is a critical point of 
$\action$ on $\Lambda$.
\end{lemma}

Now, the problem arises about which representations yield
symmetry constraints that are sufficient to imply
the existence of nontrivial (i.e., non-homographic) periodic solutions
in the equivariant loops class. As shown 
by Chenciner in \cite{chenciner}, we have to consider the problem
of collisions, coercivity and non-triviality.
We start by trying to see which conditions on $\tau$, $\rho$
and $\sigma$ might give the desired properties.

Let $\ker\tau$, $\ker \rho$ and 
$\ker \sigma$ be the kernels of the corresponding
homomorphisms 
$\rho\from G\to O(k)$, $\tau\from G \to O(2)$
and $\sigma\from G \to \Sigma_n$.
Without loss of generality we can assume that 
$\ker\tau \cap \ker\rho \cap \ker\sigma = 1$.
Moreover, assume that $g\in \ker\tau \cap \ker\rho$: 
this implies that $g\not\in \ker\sigma$.
Then, if $x(t)$ is an equivariant loop, then the restriction 
map $x(t)^g\from T^1={T^1}^g \to \Chi^g$ 
sends every point of $T^1$ to $\Chi^g$; 
since $g\in \ker \rho$ but $g\not\in \ker \sigma$, the space 
$\Chi^g$ consists entirely of collisions.
Therefore we must have $\ker\tau \cap \ker\rho = 1$ 
in order to avoid necessary collisions.

Furthermore, assume that $g\in \ker\tau \cap \ker\sigma$,
and thus $g\not\in \ker\rho$.
Again, every configuration in the orbit
$x(t)^g = x(t)$ needs to belong to $\Chi^g$, 
which is nothing but the subspace of configurations
of points in $V^g$, which is a linear proper
subspace of $V$. Thus we can consider it as a sub-problem,
and assume that $\ker\tau \cap \ker\sigma = 1$ as well.

Finally, consider $g\in \ker\rho\cap \ker\sigma$. 
Its action on $T^1$ can be a rotation or a reflection.
In case it is a rotation, we are considering $n$-bodies that 
actually tread the same loop more than once, and clearly
the problem can be solved by solving the problem concerning
loops with just one iteration. So, 
we can assume that 
every element of 
$\ker\rho\cap \ker\sigma$
acts as a reflection on $T^1$.
But if there are two distinct such elements  
$g_1\neq g_2 \in \ker\rho\cap \ker\sigma$,
then their product $g_1g_2$ would act as a rotation on $T^1$,
hence $g_1g_2$ must be trivial, i.e., $g_1=g_2^{-1} = g_2$. 
This implies that 
$\ker\rho\cap \ker\sigma$ has at most one non-trivial
element, that is it is a subgroup of order at most $2$.
If $\ker \rho\cap \ker \sigma  \neq 1$, then 
every loop in $\Lambda^G$ can be decomposed as $\gamma \gamma^{-1}$,
i.e., it is a loop that runs along a path $\gamma$ from $g(0)$ to
$\gamma(1)$ in $1/2$ of the time, 
and then from $\gamma(1)$ to $\gamma(0)$ in the second half 
of the time interval.

\begin{remark}
The equivariant loops $x\in \Lambda^G$ can be seen as 
$G/\ker \tau$-equivariant
loops $T^1 \to \Chi^{\ker \tau}$. Thus we can consider the same 
problem related to 
a finite subgroup of 
$O(2)$ 
(thus a subgroup of a dihedral group)  
and a linear subspace
$\Chi^{\ker \tau} \subset \Chi$.
Moreover, if $\Chi^{\ker \tau} \subset \Delta$ then all the loops 
are just made of collision points.
Therefore we assume that $\Chi^{\ker \tau} \not\subset \Delta$.
\end{remark}

If for a choice of $\sigma$, $\rho$ and $\tau$ 
one of the following cases occurs, we say that the action
of $G$ is \emph{degenerate}.
\begin{inparaenum}[\itshape (a)]
\item
$\ker \tau \cap \ker\sigma \cap \ker \rho \neq 1$.
\item
Every loop in $\Lambda^G$  has collisions.
\item 
There is a proper linear subspace $S$ of $E^k$ such that 
for every $t\in T^1$ and 
for every $x \in 
\Lambda^G$, the body $x_i(t)$ belongs 
to $S$.
\item 
For every loop $x(t)$ in $\Lambda^G$ 
there is a loop $y(t)$ and $k\in \ze$, $k\neq 0,\pm 1$, such that 
for every $t\in T^1$ we have $x(t)=y(kt)$. 
\end{inparaenum}

\begin{lemma}
If the action of $G$ is non-degenerate, then 
$G$ is a finite subgroup of 
$O(2)\times O(k)$ and a finite subgroup
of $O(2)\times \Sigma_n$.
\end{lemma}
\begin{proof}
Since $\ker \tau \cap \ker \rho = 1$,
the homomorphism $\tau \times \rho \from G \to 
O(2)\times O(k)$ has trivial kernel.
The same happens to the homomorphism $\tau \times \sigma$.
\end{proof}

For example, in case $k=2$, this implies that 
$G$ is a subgroup of the direct product of two dihedral groups,
and hence metabelian.
For $k=3$, $G$ is a subgroup of the direct product
of a dihedral group and a finite subgroup of $O(3)$, 
hence either $G$ is metabelian or 
it is an extension of a finite metabelian group
with with a finite subgroup of $O(3)$. 
Hence the only nonsolvable group occurs if $G$ projects onto the icosahedral
group $A_5$ in $O(3)$.

\section{Coercivity}

\begin{lemma}
\label{lemma:coercive}
The symmetric action functional $\Lambda^G$ is coercive 
if and only if $\Chi^G = 0$.
\end{lemma}
\begin{proof}
Consider the group $G/\ker \tau$. It is a finite subgroup
of $O(2)$, hence it is either a cyclic group or a dihedral group.
Let us consider first the cyclic case.  Let $c$ be a generator of 
$G/\tau$. Let $X=\Chi^{\ker\tau}$. Since $\Chi^G=X^c$, we have 
$X^c=0$. Therefore $X$ can be decomposed into irreducible
components $X=\R+\R+\dots+\R+\C+\dots+\C$,
where on the one-dimensional components $\R$ the action
of $c$ is given by $c(s)=-s$, while on the two-dimensional
components we have 
$c(z) = e^{2\pi i / l}z$ for a suitable $l\in \ze$.
Thus, using the same argument as in Bessi and Coti-Zelati  \cite{bessi},
it is possible to show that $\Lambda^G$ is coercive.
Now consider the case $G/\tau$ is dihedral. Let $h_1$ and 
$h_2$ be two generators of order $2$ of $G/\tau$.
Again, $X$ can be decomposed as 
$X=\R+\R+\dots +\R + \C \dots + \C$, where 
on the one-dimensional irreducible components
the action is either $r_1(s)=-s=r_2(s)$ or $r_1=-r_2(s)=\pm s$,
while on the two-dimensional irreducible components
$\C$ is a dihedral representation.
Thus, again exactly with same argument as \cite{bessi} 
it can be shown that there is $\alpha>0$
such that $|x|_{L^2} \leq \alpha |\dot x|_{L^2}$,
i.e., that the action functional is coercive.

For the converse, if $\Chi^G \neq 0$ let $x_0$ denote  
a loop in $\Lambda^G$ (possibly with collision) with finite action 
$\action^G$. Then $x_0 + v$, with $v\in \Chi^G$, is again a 
loop in $\Lambda^G$, with action $\action^G(x_0+v) <  \action^G(x_0)$.
But as $|v|\to \infty$ also $x_0 + v$ goes to infinity,
hence $\action^G$ is not coercive.
\end{proof}

\section{Dihedral orbits}
Consider the time circle $T^1\subset \R^2$ of radius $\frac{T}{2\pi}$,
where $T$ is the period of a periodic orbit. Let $h_1$ and 
$h_2$ be two reflections in $\R^2$ that fix two lines
forming an angle of $\frac{\pi}{l}$, with $l>1$.
Then the group $G$  generated by  $h_1$ and $h_2$ is the 
\emph{dihedral group} $D_{l}$ of order $2l$.
Consider as above a $k$-dimensional orthogonal  representation
of $G$ and an homomorphism $\sigma\from G \to \Sigma_n$
to the symmetric group of order $n!$.
This means that $h_1$ and $h_2$ act on $V=E^k$ (with a symmetry of order $2$
along a plane, a line or the origin) and on the set 
$\{1,2,\dots, n\}$ of indexes via the homomorphism $\sigma$.
The elements $\sigma(h_i)$ need to be of order $2$ in $\Sigma_n$,
whenever they are not trivial.
Given these data, $G$ acts on $\Chi$ by
\[
g (x_1,\dots, x_n) = (gx_{\sigma(1)},\dots, gx_{\sigma(n)} ),
\]
where we mean $\sigma(i) = \sigma(g)(i)$ and on $T^1$.
Then $G$ acts on the loop space $\Lambda$ by 
\[
g\cdot \gamma\from t \to g \gamma(g^{-1}t)
\]
for every $t\in T^1$ and every $\gamma\in \Lambda$.
The space $\Lambda^G$ consists of the equivariant loops. 
It is easy to see that $\Lambda^G$ is homeomorphic to the space $P$
of all the paths $\lambda\from [0,1] \to \Chi$
with the property that $\lambda(0) \in \Chi^{h_1}$ and
$\lambda(1)\in \Chi^{h_2}$. 
The homeomorphism is given by the restriction function
$r\from \Lambda^G \to P$.
The action functional can be defined in exactly
the same way on $P$, by integrating $L$ 
along $\lambda$ with a rescaled time.
Let is denote it by $\action_P$.
If $L$ is invariant with respect to the action of $G$,
then $2l \action_P r (\gamma) =  \action^G (\gamma)$,
for every $\gamma\in \Lambda^G$.
Hence $\gamma$ is a stationary point for $\action^G$ if and only if
its restriction $r(\gamma)$ is stationary for $\action_P$.
We can hence consider critical points  and local minima 
of $\action_P$ in $P$.
This is a sort of 
\emph{generalized Bolza problem}.
\begin{lemma}
Any critical point of $\action_P$  in $P$ yields 
a critical point of $\action^G$ in $\Lambda^G\subset \Lambda$, which 
is a critical point of  $\Lambda$.
\end{lemma}

\section{Minimizing on constrained paths}
More generally, assume that $h_1$ and $h_2$ are two elements of order $2$ 
acting isometrically on $E^k$ and 
$\{1,2,\dots,n\}$.
Let $X_1$ and $X_2$ be two the fixed 
subspaces  $\Chi^{h_1}$ and $\Chi^{h_2}$ of $\Chi$ 
and let $P$ denote the Sobolev space of all the paths
$\gamma\from [0,1] \to \Chi$ such that 
$\gamma(0) \in X_1$ and $\gamma(1) \in X_2$.
That is, $P = H^1( ([0,1],0,1), (\Chi,X_1,X_2) )$ is the Sobolev space 
of the $L^2$ paths $[0,1] \to \Chi$ with $L^2$ derivative 
with the constraints at the endpoints of the interval $[0,1]$.
Let $\action_P$ be defined on $P$ as above.
\begin{lemma}
Any local minimum of $\action_P$ can be extended 
to a solution (in the weak sense) 
$x\from \R \to \Chi$ which is periodic 
in a rotating frame.
\end{lemma}

Now we have to face the problem of the possible existence of collisions
in (local) minima of the action functional. 

\begin{assumption}
\label{assumption}
If $G$ is a finite group and $\Lambda^G\subset \Lambda$ is 
defined as in \ref{eq:lambdag}, then 
all local minima of the action $\action^G$ in $\Lambda^G$  are
collision-free.
\end{assumption}

We cannot term \ref{assumption} a lemma, since it has not yet been
proved in full generality. However, there is a certain evidence
that it holds for general actions.
In fact, if the group $G$ is cyclic and acts 
in the standard way on $T^1$ (that is, yielding choreographies),
then it was proved recently by Marchall (in some unpublished notes). 
The proof can be extended without significant change to some
other group actions, but will not work in full generality.
On the other hand Majer--Terracini methods on collisions singularities
\cite{serter1,mate,setate,mate95,rate,mate93,serter2,mate93a}
can be extended to the equivariant case, if $n$ is $3$ or $4$
and under some further assumptions.
But there are still some gaps in the proofs, so that we hope
to provide a complete proof in the future. 
For the purpose of this note, it suffices either to 
consider a strong--force perturbation $\Usf$ 
or to consider the numerical hint that the algorithm stops
at a non-collision loop and determines it as a global minimum.

\section{Three bodies in the plane}
Now we can start to investigate which symmetry constraints 
yield non-trivial periodic solutions.
We start with the case of $3$-bodies. 
Recall that $G$ now is a dihedral group with standard
generators $h_1$ and $h_2$.

\subsection{Three equal masses}
\label{subsec:threeeq}
We can assume that $m_i=1$ for $i=1,2,3$.
If $G$ acts on $\{1,2,3\}$ without fixed points,
then it must be $\sigma(h_1) = (1 2)$ 
and $\sigma(h_2)=(2 3)$, up to an inner permutation of the indexes.
So that the symmetric group $\Sigma_3$ is a homomorphic image of $G$.
To determine $G$ and its action on $E^2$ 
we consider now the cases for $h_1$ and $h_2$.

First case: both are rotations (of angle $\pi$) on $E^2$. Then the minimal
$G$ with this property is the dihedral group $D_3\cong \Sigma_3$ of order $6$. 
The space $\Chi^{h_1}$ is the space of all the configurations
with $x_3=0$ and $x_1 = - x_2$, while $\Chi^{h_2}$ is given by
all the configurations with $x_1=0$ and $x_2=-x_3$. 
It is clear that $\Chi^{h_1} \cap \Chi_{h_2} = 0$,
hence by \ref{lemma:coercive} the minimum $x=x(t)$ in $\Lambda^G$ exists
and is collision-free by \ref{assumption}.
Since the product $h_1h_2$ is a rotation of $T/3$ in the time circle
and acts trivially on $E^2$, we have that $x$ is a choreography.
It cannot be an Euler or Lagrange solution, hence it is a non-trivial 
choreography. 
It is possible to show that it has an ``eight'' shape.
There is the natural question, whether
it is the same as the Montgomery-Chenciner orbit or not.

Second case: $h_1$ acts on $E^2$ by reflection along a line and 
$h_2$ by rotation of angle $\pi$. 
Since the product $h_1h_2$ acts as a reflection in $E^2$ 
and as the cyclic permutation $(1 2 3)$ in $\Sigma_n$, 
$G$ needs to be the dihedral group $D_6$ of order  $12$.
The configurations in 
$\Chi^{h_1}$  are those such that $(x_1,x_2,x_3)$ is a triangle
symmetric with respect to the line fixed by $h_1$  and the configurations in
$\Chi^{h_2}$ are those such that again $x_1=0$ and $x_2=-x_3$.
This is the action described in \cite{monchen}, and the corresponding
solution 
is the figure eight choreography.

Third case: both $h_1$ and $h_2$ act on $E^2$ by a reflection along 
a line ($l_1$ and $l_2$ respectively). 
If $l_1=l_2$, then the product $h_1h_2$ acts trivially
on $E^2$, hence the minimal $G$  is the dihedral group $D_6$.
The minimum exists and numerical experiments 
let one guess that it is the Lagrange orbit.
Otherwise, $l_1$ and $l_2$ intersect with an angle  $\pi/q$,
with $q>1$ integer.
The minimal $G$ is therefore $D_q$ if $q$ is $0 \mod 3$
and $D_{3q}$ otherwise. Since the Lagrange orbit
belongs to $\Lambda^G$ it can be the minimum.
We did not check whether it is a minimum for every $q$ or not.

\subsection{Two equal masses}
\label{subsec:3_2eq}
Now assume that the first two masses are equal ($m_1=m_2$).
We again list the possible cases.
Without loss of generality we can assume 
$\sigma(h_2) = (1 2)$, since at least one of $h_1$ and $h_2$
needs to  act non-trivially on the indexes.

First case: $\sigma(h_1) = (1 2)$. 
If both $h_1$ and $h_2$ act rotating on $E^2$, 
then $\Chi^{h_1} = \Chi^{h_2} = \Chi^{G}$, and the functional
is not coercive.
If $h_1$ acts by reflection and $h_2$ by rotation, then 
again the functional is not coercive, 
and the same is true if they act by reflecting along the 
same line.
So it is left to check the case in which they act 
on $E^2$ by reflection along two distinct lines.
In this case the functional is coercive, and $\Lambda^G$ contains
the Lagrange solution so that it is of no interest.

Second case: $\sigma(h_1) = ()$.
Since $h_1$ does not move the indexes, 
to avoid collisions it is  is necessary that 
$h_1$ does not act on $E^2$  with a rotation.
It cannot have a trivial action, since otherwise $\Chi^{h_1} = \Chi$
and the orbit would not be dihedral, so that it will be 
a reflection. Now, it is left to determine the action of 
$h_2$. If $h_2$ is a rotation of angle $\pi$, then 
the action functional is not  coercive.
\begin{remark}
We can restrict the space of paths considering only paths
$x$ 
with a prescribed order of the configuration $x(0)$.
If we look at the configurations
such that $x_3(0)$ does not lie between $x_1$ and $x_2$,
then $\action$ is coercive.  
In a strong-force settings a collision-free minimum need to exist.
Numerical experiments show that 
such minima might exist even for a potential of type $1/r^{a}$, with
$a\geq 1.3$ (see figure \ref{fig:1}).
\ref{fig:1}.
\begin{figure}
\begin{center}
\includegraphics[width=4truecm]{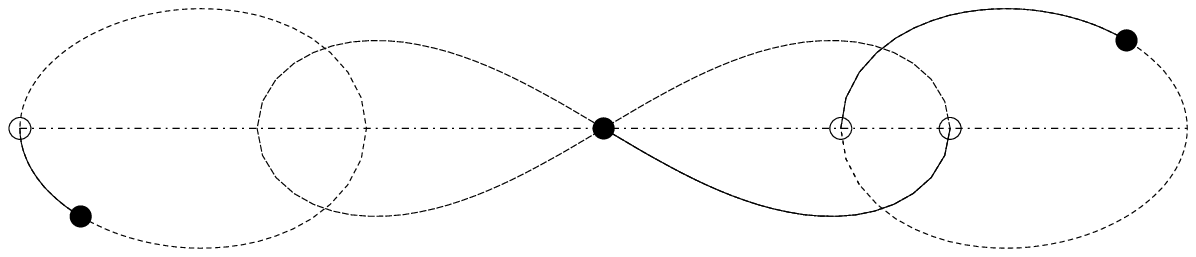}
\caption{Two equal masses }
\label{fig:1}
\end{center}
\end{figure}
This is the solution of braid $b_1^2b_2^{-2}$ of Moore \cite{moore},
pag. 3677.
\end{remark}

\begin{figure}
\begin{center}
\includegraphics[width=4truecm]{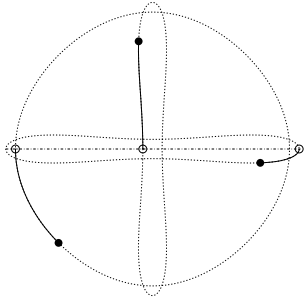}
\caption{$m_1=m_2=m_3=1$}
\label{fig:3b}
\end{center}
\end{figure}

Otherwise, $h_2$ is a reflection along a line.
If  this line coincides with the line fixed by
$h_1$, then again the functional is not coercive,
since there are orbits in which $x_1$ and $x_2$ rotate in a circle
very far from $x_3$.
So we can assume that the two lines intersect with an angle
$0<\alpha \leq  \pi/2$. 
If the angle $\alpha$ is equal to $\pi/2$ then again $\action$ is not 
coercive, so we assume  $0<\alpha<\pi/2$.
Now the functional is coercive, and there is a minimum.
If $q$ is an integer, then 
the minimal group $G$ is $D_q$ if $q$ is even
and $D_{2q}$ if $q$ is odd.
At $t=2T/q$ the configuration is the same as the 
configuration 
at $t=0$ with the two bodies interchanged and rotated by
an angle of $2\pi/q$; at $t=4T/q$ it 
is exactly the configuration at $t=0$ rotated of an angle
$4\pi/q$.
If $q>2$ is not an integer, then one obtains an orbit
periodic with respect to a rotating frame.
Unfortunately the Euler orbits belong to this class,
so that the minimum can be achieved on a Euler solution.
Some numerical simulations give a hint that this is not the case:
orbits like the one in figure \ref{fig:3b} can be found
with constrained optimization techniques, with an action 
less than the action of the corresponding Euler orbit.

To prove the existence of such orbits, provided that because 
of \ref{assumption} there are no collisions in a minimizing
orbit, it suffices to use the following level estimates.
We are going to compare the levels of the action
of suitable symmetric orbits   with the action of the Euler orbit.
Let $m_1=m_2 = 1$, $m_3=m>0$ be the masses.
Let $c$, $r_0$ and $l$ be three positive constants (to be determined),
and consider the path in $P$ determined by the equations
\begin{equation}
\label{orbit1}
\begin{split}
x_1(t) & = 
le^{i\theta t} + (r_0+ct)e^{i(\theta - \frac{\pi}{2})t} 
\\ 
x_2(t) & =  
le^{i\theta t} - (r_0+ct)e^{i(\theta - \frac{\pi}{2})t} 
\\
x_3(t) & = 
-\frac{2l}{m} e^{i\theta t}.
\\ 
\end{split}
\end{equation}

The kinetic contribution of $(1)$ and $(2)$ to the action $\action(x)$  is
\[
\begin{split}
K_1 + K_2 = 
1/12\,{c}^{2}{\pi }^{2}+1/4\,{{r_0}}^{2}{\pi }^{2}+  {r_0}\,{
\theta}^{2}c 
- c{r_0}\,\theta\,\pi +{c}^{2}+1/3\,{c}^{2}{\theta}^{2
}-{{r_0}}^{2}\theta\,\pi +\\
 + {l}^{2}{\theta}^{2}-1/3\,{c}^{2}\theta
\,\pi +{{r_0}}^{2}{\theta}^{2}+1/4\,{r_0}\,{\pi }^{2}c. \\
\end{split} 
\]
The kinetic term coming from $(3)$  is  simply
\[
K_3 = 2\frac{l^2 \theta}{m}.
\]
Now consider the terms corresponding to the potential.
The term corresponding to the interaction between $(2)$ and $(3)$ 
is equal to
\[
U_3 = \frac{1}{2c} \log(1+\frac{c}{r_0}).
\]
Now, the term of the interaction between $(1)$ and $(3)$ is bounded 
by
\[
U_2 \leq 
m \left( (r_0 + c)^2 + l^2(1+2/m)^2 \right)^{-1/2},
\]
and a similar inequality holds for the term $(2)-(3)$:
\[
U_1 \leq 
\dfrac{m^2}{m(l-r_0)+2l}.
\]
Let $\action_D = K_1+K_2+K_3 + U_1 + U_2 + U_3$ denote the Lagrangian action of the path
\ref{orbit1}.
The action of the Euler solution with the body $x_3$ in the center of mass
is 
\[
\action_E =
\frac{3}{2} \left[ 
(1/2 + 2m)^2(\pi/2 -\theta)
\right]^{1/3}.
\]
With some computations it is possible to simplify the difference as 
\begin{equation}
\label{eq:difference}
\begin{split}
\action_D-\action_E =
l^2 \theta^2 (1+2/m) + c^2 
+ \frac{1}{12} (3r_0^2 + c^2 + 3cr_0)(\pi-2\theta)^2
+ \dfrac{m^2}{m(l-r_0)+2l} +
\\
+ \dfrac{m}{ \sqrt{ (r_0 + c)^2 + l^2(1+2/m)^2 }} +
\frac{1}{2c} \log(1+\frac{c}{r_0}) 
 -
\frac{3}{2} \left[ 
(1/2 + 2m)^2(\pi/2 -\theta)
\right]^{1/3}.  \\
\end{split}
\end{equation}

Now, let $D\subset \R^2$ the domain of all the pairs $(m,\theta)$
such that 
\begin{equation}
\label{eq:defD}
\begin{split}
\mathrm{inf} \{ \action_D - \action_E \st l>r_0>0, l>c>0 \} < 0.
\end{split}
\end{equation}

The following proposition is a trivial consequence of the definition
of $D$.
\begin{propo}\label{propo:first}
If $(m,\theta) \in D$ then 
there are $l$,$r_0$ and $d$ such that 
the action of  the orbit \ref{orbit1}
is less than the action of the Euler orbit. Therefore 
non-homographic dihedral orbits exists for every $(m,\theta)\in D$.
\end{propo}

\begin{propo}\label{propo:last}
The set $D$ is a non-empty open subset of $\R^2$.
\end{propo}
\begin{proof}
Since $\inf$ is upper semi continuous, $D$ is open. We only need to 
show that it is not empty: 
let $m=2$, $\theta=\pi/8$,
$l=1$, $r_0=0.4$, $c=0.3$. 
Evaluating for such values an approximation of \ref{eq:difference} yields
$-.124390105 < 0$.
A candidate for the corresponding  minimum can be seen
in figure \ref{fig:candidate}.
\begin{figure}
\begin{center}
\includegraphics[width=4truecm]{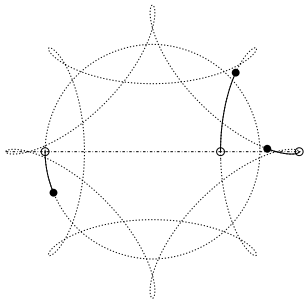}
\caption{$m_1=m_2=1$, $m=m_3=2$ and $\theta=\pi/8$.}
\label{fig:candidate}
\end{center}
\end{figure}
\end{proof}

\begin{remark}
Of course propositions \ref{propo:first} and  \ref{propo:last}
do not imply that the orbit in figure \ref{fig:3b} exists.
The minimization in that case has been done with a piecewise linear
path, which at the moment cannot be reproduced symbolically,
whose action is less than $\action_E$. The orbit in figure
\ref{fig:3b} is the same as the orbit
with braid $b_1^2b_2b_1^{-2}b_2$, found numerically
by Moore in \cite{moore}.
\end{remark}

\section{Four bodies in the plane}
\label{sec:4}
Now we analyze in the same way the situation of $4$ bodies in the plane.
Let again $h_1$ and $h_2$ be the reflections in $T^1$ generating $G$.

\subsection{Four equal masses}
\label{subsec:4eq}
Assume that $m_i=1$, $i=1,\dots, 4$.
We want that $\sigma(h_1)$ and $\sigma(h_2)$ generate
a subgroup of $\Sigma_4$ that acts transitively
on the indexes $\{1,2,3,4\}$.

First case: $\sigma(h_1) = (1 2)$. Then the only transitive
subgroup up to inner automorphism of $\Sigma_4$ is
given by the choice $\sigma(h_2) = (1 3)(2 4)$.
The action of $h_1$ on the plane $E^2$ 
cannot be a rotation (since fixes the indexes $3$ and $4$).
Since it cannot be trivial (otherwise 
at $t=0$ a collision is not avoidable) 
it needs to be a reflection along a line $l_1$.
Now, if $h_2$ acts by reflection along a line $l_2$, there
are the following sub-cases.
If $l_2=l_1$, then the homographic solution of a rotating
square can be a minimum (actually, apparently it is the minimum),
so we do not consider this case.
If $l_2$ and $l_1$ meet at an angle $\pi/4$, then 
the action is not coercive (there is a big square with stationary
masses). On the other hand even 
if the angle is different from $\pi/4$ then 
the homographic solution of the rotating square  can be a minimum,
hence we do not consider this case too (even if it might be possible
that the minimum is achieved by a non-homographic orbit). 

So, it is only left the case in which 
$h_2$ acts by rotation of angle $\pi$ on $E^2$.
Since $\Chi^{h_1} \cap \Chi^{h_2} \neq  0$, 
by lemma \ref{lemma:coercive} the action is not coercive.

Now consider the second case  $\sigma(h_1)=(1 2)(3 4)$. 
Then the only choice of $\sigma(h_2)$ that has has not yet been
considered for a transitive action is $\sigma(h_2) = (1 3)(2 4)$.
This time $h_1$ can act on $E^2$ either as a reflection  or as a rotation,
and the same holds for $h_2$.
If one is a reflection and one is a rotation,
then 
the action is not coercive, since a sequence of increasing stationary
squares can have as small as possible action.
If both act as rotations, then the group $G$ is equal 
to the group generated by $\sigma(h_1)$ and $\sigma(h_2)$,
i.e., the elementary abelian group $\ze_2^2$ of order $4$.
The 
resulting symmetric orbit can be two coupled Kepler orbits,
hence there is no coercivity.
So it is left the case in which both are reflections along lines
$l_1$ and $l_2$ in $E^2$.
If the lines $l_1=l_2$ coincide, 
then again  it is possible to see that the functional not is coercive,
by taking two symmetric Keplerian orbits that have increasing
distance (i.e., by applying lemma \ref{lemma:coercive}).

\begin{figure}
\begin{center}
\includegraphics[width=4truecm]{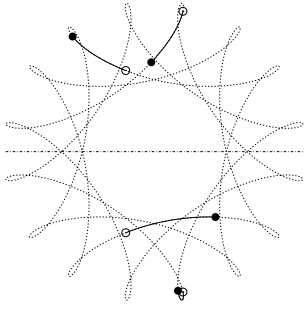}
\caption{Four equal masses with $l_1/l_2$-angle $\pi/8$.}
\label{fig:nonhomographic}
\end{center}
\end{figure}

Otherwise, if they are orthogonal, then again it is easy to see that 
the functional is not 
coercive. If they are not orthogonal  then
the functional is coercive, but in the class of symmetric paths there are 
the homographic orbits of rotating squares.
Of course the question arises whether the homographic
orbits achieve the  minimum or not.
Numerical simulations lead to  think that 
the minimum might be achieved by non-homographic orbits,
like the one depicted in figure \ref{fig:nonhomographic}.
It could be of some interest to prove some estimates
like those in the previous section, to actually prove
or disprove their existence. This is true also for other 
examples listed below, and we will not rise the question again.

\subsection{Three equal masses}
Assume now that there are three equal masses $m_1=m_2=m_3=1$ 
and a fourth mass $m_4=m$. Then the subgroup
of $\Sigma_4$ generated by $\sigma(h_1)$ and $\sigma(h_2)$
needs to act transitively only on the set $\{1,2,3\}$.
The only possibility, up to rearranging the indexes,
is $\sigma(h_1)=(1 2)$ and $\sigma(h_2) = (1 3)$,
like in the case of $3$ bodies.
Since both fix two indexes, the actions of $h_1$ and 
$h_2$ on $E^2$ need to be  reflections along the lines
$l_1$ and $l_2$ respectively.
If the lines coincide, then a rotating triangle with the mass $(4)$
in the center can be the minimum. 
If the angle is $\pi/3$, then the problem is no longer coercive,
since any constant equilateral triangle with bodies $(1)$, $(2)$ and $(3)$
with $(4)$ in the center is symmetric with respect to this action.
On the other hand, such a triangle when rotating at a suitable speed 
always belongs to the set of 
symmetric loops $\Lambda^G$, however the two lines intersect.
We do not know whether it is a minimum in $\Lambda^G$.

\subsection{Two pairs of equal masses}
\label{subsec:twopairs}
Assume that $m_1=m_2=1$ and $m_3=m_4=m$.
Since at least one from $\sigma(h_1)$  and $\sigma(h_2)$
is non-trivial, we can assume that 
$\sigma(h_1) = (1 2)$  or $\sigma(h_1) = (1 2)(3 4)$.

In the first case, $\sigma(h_1)=(1 2)$, necessarily it must be
$\sigma(h_2) =(3 4)$.
The only possible action of $h_1$ and $h_2$ on $E^2$ 
is given by reflections along lines $l_1$ and $l_2$.
If the lines coincide, then the functional is not coercive
by \ref{lemma:coercive}.
It is not
coercive also if they are orthogonal: a square in increasing
size can give a sequence going to infinity with bounded
action.
If the lines meet with an angle $\pi/q$, then it is coercive,
but a rotating central configuration with the masses 
at the vertexes of a parallelogram belongs to $\Lambda^G$,
and hence it can be the homographic minimum.
Again, as above, the question arises whether the minimum
is homographic or not.

\begin{figure}
\begin{center}
\includegraphics[width=4truecm]{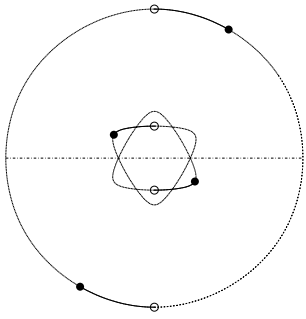}
\caption{$m_1=m_2=1$ and $m_3=m_4=4$, with $D_3$-symmetry.}
\label{fig:rough}
\end{center}
\end{figure}

Consider now the second case, $\sigma(h_1) = (1 2)(3 4)$.
There are three possibilities for $\sigma(h_2)$: 
up to rearranging indexes,
the trivial $()$, or $(1 2)$ or $(1 2)(3 4)$.
Consider $\sigma(h_2) = ()$. Then $h_2$ must act on 
$E^2$ by reflection along a line $l_2$.
However $h_1$ acts on $E^2$, as a rotation or as a reflection,
it is possible to find 
a rotating collinear central configuration belonging to
$\Lambda^G$.
So we consider the next case, $\sigma(h_2) = (1 2)$.
Again $h_2$ needs to act as a reflection along a line $l_2$,
and rotating collinear configurations now cannot belong to $\Lambda^G$.
If $h_1$ acts by rotation, then a rotating parallelogram belongs
to $\Lambda^G$, hence we consider only the case of $h_1$ acting
by reflection along a line $l_1$.
If the two lines coincide or are orthogonal, then the action functional is not 
coercive. Otherwise it is coercive, and hence
there is a minimum, which is collisionless due to theorem
\ref{assumption}.
Can it be homographic? No: at the time $t=0$ (i.e., the time in $T^1$
fixed by $h_1$) the lines through $(1)\mbox{--}(2)$ and
$(3)\mbox{--}(4)$ are parallel (both orthogonal to $l_1$),
while at time $1$ (i.e., the time in $T^1$ fixed by $h_2$)
they are orthogonal.
Such orbits can be described as follows: two masses in a 
roughly Keplerian orbit outside, and two masses in a 
retrograde approximate Keplerian orbit inside,
like in figure
\ref{fig:rough}.

Now it is left the case $\sigma(h_1) = \sigma(h_2) = (1 2)(3 4)$.
Since we are assuming the action of $G$ on $T^1$ to be faithfully 
dihedral, $h_1$ and $h_2$ must act on $E^2$ in different ways,
so that at least one of them acts as a reflection.
Let us assume that $h_1$ acts by reflecting along a line $l_1$.
If $h_2$ acts by rotation, then the functional is not coercive;
if it acts as a reflection along a line $l_2\neq l_1$, 
then homographic orbits belong to $\Lambda^G$,
so that this case is of minor interest.

\section{Orbits with an additional central symmetry}
Consider the symmetries in the previous sections.
If we can find an element $\sigma_3$ in $\Sigma_n$ of order $2$
that fixes at most one index and commutes
with $\sigma(h_1)$ and $\sigma(h_2)$,
we can consider the following additional central symmetry:
the symmetry group is $G\times \ze_2$ (where $G$ is the group
in the example under consideration), where the direct
factor $\ze_2$ is generated by $h_3$; this element
acts trivially on $T^1$, acts as a rotation of angle $\pi$
in the plane $E^2$, and is sent to $\sigma_3$ by
the homomorphism $\sigma$.
This means that for every $t\in T^1$ the configuration
at time $t$ is in $\Chi^{h_3}$, that is, bodies
(with the same mass) in the same cycle in $\sigma_3$ 
are symmetric with respect $0\in E^2$, while the possible
body with index fixed by $\sigma_3$ lies in $0\in E^2$.
If $n=3$, then such orbits are trivial, since 
they need to be always collinear.
So consider the case $n=4$. We can analyze the cases
exploited in section \ref{sec:4}
to see when these conditions are fulfilled,
and if the additional central symmetry $h_3$ yields
non-homographic orbits.
We omit the details of this case-by-case analysis, and 
exhibit only a particular family of dihedral orbits.

\subsection{Four equal masses}
\label{subsec:4eqp}
In section \ref{subsec:4eq} 
consider the case of 
$\sigma(h_1) = (1 2)(3 4)$ 
and $\sigma(h_2) = (1 3)(2 4)$,
where $h_1$ and $h_2$ act on $E^2$ by reflection
along different lines $l_1$ and $l_2$.
If $l_1$ and $l_2$ are orthogonal, then 
the functional is not coercive, even 
when we add the additional symmetry $\sigma_3 = (1 2)(3 4)$.

\begin{figure}
\begin{center}
\includegraphics[width=4truecm]{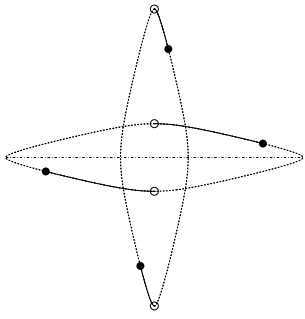}
\caption{$m_i=1$, with $(D_4\times \ze_2)$-symmetry}
\label{fig:mac}
\end{center}
\end{figure}

\begin{figure}
\begin{center}
\includegraphics[width=4truecm]{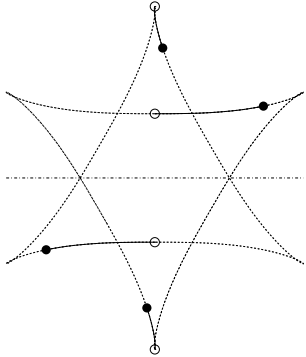}
\caption{$m_i=1$ with $(D_6\times \ze_2)$-symmetry.}
\label{fig:caseq6}
\end{center}
\end{figure}

In case $l_1$ and $l_2$ meet at an angle $\pi/q$, with 
$q>2$, we can avoid the homographic solution
again by the same additional symmetry $\sigma_3$. 
The group $G$ acts faithfully on $T^1$  and 
is equal to the dihedral group of order $2k$, 
where $k$ is the least common multiple of $2$ and $q$.
Thus we obtain an infinite family of 
periodic orbits in the $4$-body problem with equal masses.
We can see the case 
$q=4$ in figure
\ref{fig:mac}
and  
$q=3$ in figure 
\ref{fig:caseq6}.
The orbit of figure \ref{fig:mac} is very likely the orbit  found 
by Chen \cite{chen}.

\subsection{Two pairs of equal masses}
\label{subsec:22eqp}
\begin{figure}
\begin{center}
\includegraphics[width=4truecm]{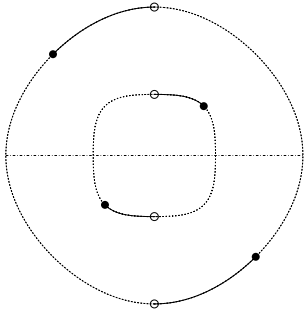}
\caption{$m_1=m_2=1$, $m_3=m_4=2$, with $(D_4\times \ze_2)$-symmetry.}
\label{fig:arch}
\end{center}

\end{figure}

The periodic solutions of section
\ref{subsec:twopairs}
can be endowed with the additional symmetry given by
$\sigma_3=(1 2)(3 4)$. So we consider 
$\sigma(h_1) = (1 2)(3 4)$,
$\sigma(h_2) = (1 2)$,
the action of $h_1$ and $h_2$ on $E^2$ is by reflection
along two lines $l_1$ and $l_2$ that intersect
at an angle $\pi/q$ with $q>2$.
We can see the case $m_1=m_2=1$, $m_3=m_4=2$ 
and $q=4$ in figure \ref{fig:arch}.

\section{Some plane choreographies for $n>3$ bodies}
\label{sec:choreo}
As shown in \cite{chenciner},
it is not difficult to generalize the eight-shaped choreography
of Montgomery-Chenciner to the case of $n>3$ odd bodies with 
equal masses.
Consider 
the following permutations on $\{1,2,\dots,n\}$:
\begin{equation}
\begin{split}
\sigma_1\from i \to n-i \mod n\\
\sigma_2\from i \to n-i+1 \mod n,\\
\end{split}
\end{equation}
where we understand that $0 \equiv n \mod n$.
That is, for $n=3$ we have 
$\sigma_1 = (1 2)$, $\sigma_2=(1 3)$;
for $n=5$ we have 
$\sigma_1=(1 4)(2 3)$ and 
$\sigma_2=(1 5)(2 4)$;
for $n=7$ it is 
$\sigma_1=(1 6)(2 5)(3 4)$
and $\sigma_2=(1 7)(2 6)(3 5)$.
The product $\sigma_1 \sigma_2$ ($\sigma_2 \sigma_1$ in functional
notation) sends $i$ to $i+1 \mod n$, i.e., 
$\sigma_1 \sigma_2$ is the cyclic permutation
$(1 2 \dots n)$.
Thus the subgroup generated by $\sigma_1$ and $\sigma_2$
is a dihedral group of order $2n$. We can define
$\sigma_1$ and $\sigma_2$ in a geometrical way:
consider a regular $n$-gon with consecutive vertices
$(1),(2), \dots, (n)$. Then $\sigma_1$ is the reflection
with axis through the vertex $(n)$ and $\sigma_2$ the 
reflection fixing the vertex $(i)$ with  $i=(n+1)/2$.
We will of course choose $\sigma(h_i) = \sigma_i$,
with $i=1,2$, where $h_i$ are the generators
of the symmetry group $G$ as above (it is not assumed 
that the homomorphism $\sigma$ is a monomorphism;
it will depend on the choice of the action of $h_1$ add $h_2$
on $E^2$).

\begin{figure}
\begin{center}
\includegraphics[width=4truecm]{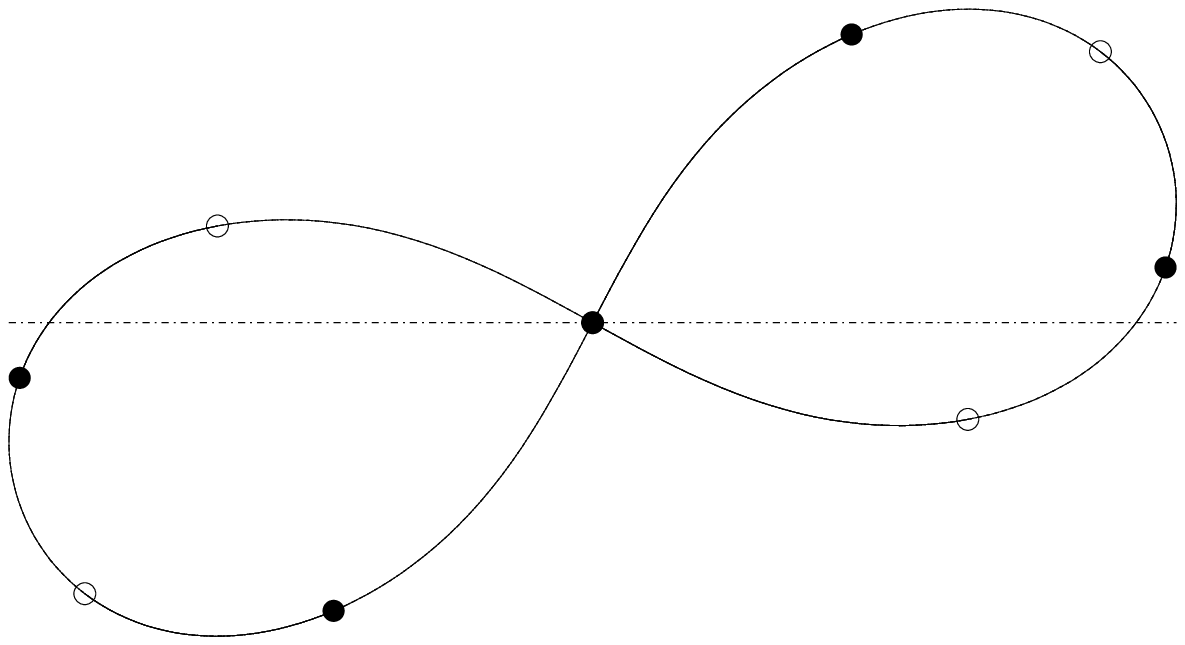}
\caption{$5$-body choreography}
\label{fig:chor5}
\end{center}
\end{figure}

Hence, it is only left 
to choose the action of $h_1$ and $h_2$ 
on the plane $E^2$.
Again, $h_1$ and $h_2$ need to be of order two,
hence they can be either rotations of angle $\pi$
or reflections along lines $l_1$ or $l_2$.
Let us consider first the case in which $h_1$ and $h_2$
are both the rotation of angle $\pi$.
The group $G$ is therefore the dihedral group 
of order $2n$.
By 
lemma \ref{lemma:coercive}
the functional is coercive, hence
it attains the  minimum;
by \ref{assumption},
the minimum is collision-free.
It is only left to show that this minimum
is not homographic.
This is easy, since at time $t=0$
the body $(n)$ is in the origin $0\in E^2$,
while at time $T/(2n)$ 
the body in the origin is the body $(i)$ with $i=(n+1)/2$
(figure  \ref{fig:chor5}).

\begin{figure}
\begin{center}
\includegraphics[width=4truecm]{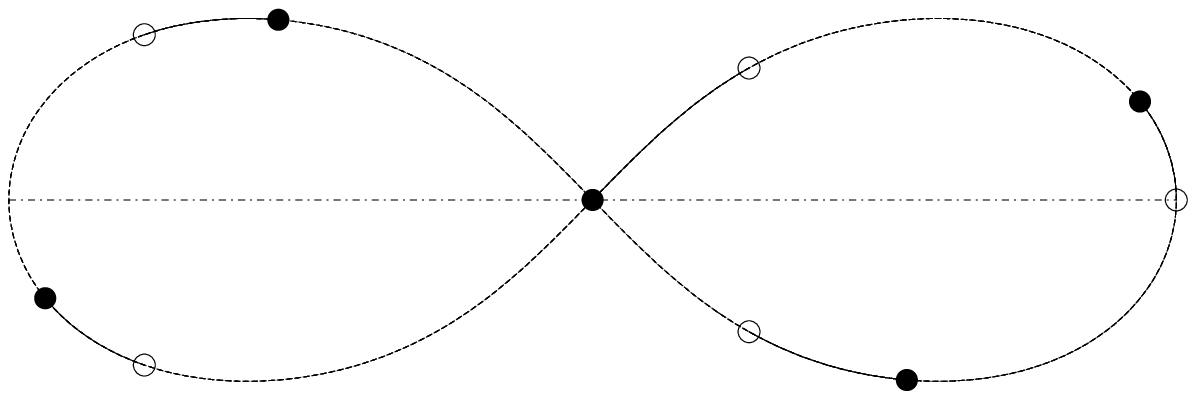}
\caption{$5$-body choreography}
\label{fig:chor5bis}
\end{center}
\end{figure}

If $h_1$ acts by rotating and $h_2$ by a reflection,
then the same results hold,
only this time the group $G$ is the dihedral group
of order $4n$.
Again, we have a choreography with $n$ equal masses. 
It might be interesting to see whether 
it is the same  as the choreography with $h_i$ rotations or not
(we have tried some simulations, obtaining something like 
figure  \ref{fig:chor5bis}, which is very similar to 
figure \ref{fig:chor5}).

In case $h_1$ and $h_2$ both act 
by reflections, then clearly a rotating regular $n$-gon
belongs to $\Lambda^G$ (if the functional is coercive),
so that the minimum can be homographic and we cannot
apply the previous results.


\section{Remarks}
\begin{remark}
The methods used in the paper for proving the existence
of non-homographic periodic orbits are quite simple, 
once the collision assumption is proved, 
and can be extended in a straightforward way to 
the case of $n>4$ bodies and non--dihedral groups.
A full classification of non-degenerate actions of dihedral groups
or abelian extensions of dihedral groups is not difficult,
and will be the content of a forthcoming paper.
\end{remark}

\begin{remark}
So far, in this note we have considered explicitly 
only planar periodic orbits. This is not a serious restriction,
since the only change needed to deal with the case $E^3$ 
of non-planar periodic orbit is to add a one-dimensional irreducible
representation of $G$ to the representation $\rho$ (there are 
no $3$-dimensional irreducible representation of the dihedral group).
The known periodic orbits in the space (such as 
the Chenciner--Venturelli ``hip--hop'' orbit \cite{chenven})
obtained by symmetry constraints use a $3$-dimensional representation
which can be decomposed into an irreducible $2$-dimensional
representation and a $1$-dimensional representation.

In fact, the method is quite simple. Consider a periodic orbit
arising from the method given above, with generators $h_1$ 
and $h_2$ of $G$. Then $h_1$ and $h_2$ can act on $E^2\subset E^3$
by this action, and on the third
orthonormal coordinate of $E^3$ as $\pm 1$.
If they act both trivially, then the problem will not be coercive. 
Otherwise, we can have three choices ($(+,-)$,$(-,+)$ and $(-,-)$)
that will yield non-degenerate coercive actions on $\Chi$.
In some cases $(+,-)$ and $(-,+)$ will yield the same 
periodic orbit, up to a time shift, but in general we will
get three periodic orbits in the space.
Numerical simulations can be done exactly as in the planar case:
we obtained some interesting analogues of the ``hip--hop''
orbit.
\end{remark}
\begin{remark}
The program used for the simulations is rather simple. 
We consider a PL discretization of the loops, and 
so we obtain a finite dimensional space $\Lambda$.
Then by relaxation dynamics on $\Lambda^G$ (that is,
a very simple gradient method) and a 
random method for avoiding poor progress 
(that is, we restart the relaxation process after a
random small variation within $\Lambda^G$ if the progress is not good,
until it is apparent that the program is in a local minimum)
we obtain an approximation of the minimum.
Now we can compare the action functional on such a path (that can be computed
with a reasonable precision), and compare with other 
known values (like homographic solutions).
The language used was FORTRAN 95 with double precision arithmetic, 
and the NETLIB SLATEC
library for the ODE solver and error-handling routines.
The figures are produced with GNUPLOT run on the raw data files.
\end{remark}
\begin{remark}
This is a very preliminary report. Not only the collision
assumption \ref{assumption} is still present, but also the program used
is quite bad designed and 
has a very poor performance. 
While writing the program I was more
concerned about flexibility, robustness and simplicity, than performance
or very good approximation of the solutions.
Thus the algorithm is very slow 
and does not give a very good approximation of the orbits.
Moreover, I did not compute the linear stability of the orbits,
nor I used the more efficient approaches of variational
optimization techniques in symmetric periodic problems,
available in the literature. 
\end{remark}


\def\cfudot#1{\ifmmode\setbox7\hbox{$\accent"5E#1$}\else
  \setbox7\hbox{\accent"5E#1}\penalty 10000\relax\fi\raise 1\ht7
  \hbox{\raise.1ex\hbox to 1\wd7{\hss.\hss}}\penalty 10000 \hskip-1\wd7\penalty
  10000\box7} \def\cprime{$'$} \def\cprime{$'$} \def\cprime{$'$}
  \def\cprime{$'$}


\end{document}